\documentclass{amsart}
\usepackage{xspace,amssymb,amsfonts,euscript,ifpdf}
\ifpdf
\usepackage{pdfsync}
\fi
\author{Ben Webster}
\title[Cramped subgroups]
  {Cramped subgroups and generalized Harish-Chandra modules}
\address{Department of Mathematics\\
         University of California, Berkeley\\
         Berkeley, CA 94720}
\curraddr{Department of Mathematics\\ Massachusetts Institute of Technology\\ 77 Massachusetts Avenue\\ Cambridge, MA 02139}
 \email{bwebster@math.mit.edu}
 \urladdr{http://math.mit.edu/\~{}bwebster}
\thanks{This material is based upon
  work supported under a National Science Foundation Graduate Research 
  Fellowship and partially supported by the RTG grant DMS-0354321.}
\begin{document}
\newtheorem{thm}{Theorem}[section]
  \newtheorem*{thm*}{Theorem}
  \newtheorem{prop}[thm]{Proposition}
  \newtheorem{lem}[thm]{Lemma}
  \newtheorem{cor}[thm]{Corollary}
  \newtheorem{conj}{Conjecture}
 
 \theoremstyle{remark}
 \newtheorem{defi}{Definition}

 \newcommand{\nc}{\newcommand}
 \newcommand{\renc}{\renewcommand}
 \nc{\Ann}{\mathrm{Ann}}
 \nc{\Rmodgr}{R-\mathsf{modgr}}
 \nc{\ghc}{generalized Harish-Chandra\xspace}
 \nc{\ghcm}{\ghc module\xspace}
 \nc{\Hdim}[1]{\EuScript{D}_{#1}}
 \nc{\irr}{\mathrm{vag}\,}
 \nc{\Mdis}[1]{\EuScript{N}_{#1}}
 \nc{\supp}{\mathrm{supp}}
 \nc{\asupp}{\mathrm{asupp}}
 \nc{\wlat}{\mc X}
 \nc{\ga}{\gamma}
 \nc{\Si}{\Sigma}
 \nc{\al}{\alpha}
 \nc{\C}{\mathbb{C}}
 \nc{\de}{\delta}
 \nc{\fr}[1]{\mathfrak{#1}}
 \nc{\Hom}[3]{\mathrm{Hom}_{#3}(#1,#2)}
 \nc{\im}{\mathrm{im}\,}
 \nc{\la}{\lambda}
 \nc{\mc}[1]{\mathcal{#1}}
 \renewcommand{\O}{\mathcal{O}}
 \nc{\om}{\omega}
 \nc{\R}{\mathbb{R}}
 \nc{\Stab}{\mathrm{Stab}}
 \nc{\Z}{\mathbb{Z}}
 \nc{\comG}{K}
 \nc{\comH}{L}
 \nc{\comg}{\fr k}
 \nc{\comh}{\fr l}
 \nc{\ssG}{G}
 \nc{\ssH}{H}
 \nc{\ssg}{\fr g}
 \nc{\ssh}{\fr h}
 \nc{\tg}{\tilde{\ssg}}
 \nc{\del}{\delta}

\theoremstyle{remark}
\newtheorem*{question}{Question}

\begin{abstract}
  Let $\ssG$ be a reductive complex Lie group with Lie algebra $\ssg$.
  We call a subgroup $\ssH\subset \ssG$ {\bf cramped} if there is an
  integer $b(\ssG,\ssH)$ such that each finite dimensional
  representation of $\ssG$ has a non-trivial invariant subspace of
  dimension less than $b(\ssG,\ssH)$. We show that a subgroup is cramped
  if and only if the moment map $T^*(\comG/\comH)\to\comg^*$ is
  surjective, where $\comG$ and $\comH$ are compact forms of $\ssG$ and
  $\ssH$. We will use this in conjunction with sufficient conditions for
  crampedness given by Willenbring and Zuckerman \cite{WZ04} to prove a
  geometric proposition on the intersections between adjoint orbits and
  Killing orthogonals to subgroups.
  
  We will also discuss applications of the techniques of symplectic
  geometry to the generalized Harish-Chandra modules introduced by
  Penkov and Zuckerman \cite{PZ04}, of which our results on crampedness
  are special cases.
\end{abstract}

\maketitle

\section{Introduction}
\label{sec:introduction}

Throughout this paper, we fix an inclusion of semi-simple complex Lie
groups $\ssH\subset \ssG$.  We will assume throughout that $G$ and $H$
are connected.  Further, we fix Cartan and Borel subgroups $T\subset
B\subset G$, such that the intersection $B\cap \ssH$ is a Borel subgroup for
$\ssH$.  Let $\comG,\comH$ be the fixed points, on $\ssG$ and $\ssH$
respectively, of the Cartan involution associated to $T\subset B$.
These are maximal compact subgroups of $G$ and $H$, respectively, and
$S=T\cap\comG$ is a maximal torus of $\comG$.

We let $\wlat_+^\ssG$ and $\wlat_+^\ssH$ be the cone of dominant
integral weights of the respective groups.  For
$\la\in\wlat_+^G,\mu\in\wlat_+^H$, let $V_\la$ and $W_\mu$ be the
corresponding finite-dimensional irreducible representations.

For each weight $\la\in\wlat_+^G$, we can consider the representation
$V_\la$ restricted to $\ssH$.  Let $b(\la)$ be the minimum among the
dimensions of non-trivial $H$-invariant subspaces of $V_\la$, that is,
the dimension of the smallest irreducible $\ssH$-representation which
appears as a subspace of $V_\la$.  Then we have an invariant of the
inclusion $H\hookrightarrow G$ given by
\begin{equation*}
b(\ssG,\ssH)=\sup\{b(\la)\vert\la\in\wlat_+^G\}. 
\end{equation*}
That is, $b(\ssG,\ssh)$ is the smallest integer such that every non-trivial $G$-representation has a non-trivial $H$-invariant subspaces $\leq b(\ssG,\ssh)$

\begin{defi}
  We call the subgroup $\ssH$ or corresponding
  inclusion {\bf cramped} if $b(\ssG,\ssH)$ is finite.
\end{defi}
An inclusion of Lie algebras $\ssh\subset\ssg$ will also be called
{\bf cramped} if the corresponding inclusion of groups
$\ssH\subset\ssG$ where $\ssG$ is simply-connected is cramped.

Note that the notion of crampedness requires a choice of a category of
representations. By convention, we will consider the category of
finite-dimensional complex algebraic representations if the group in
question is complex reductive, and the category of finite-dimensional
smooth complex representations if the group is compact.  Since
restriction from a complex reductive group to its maximal compact,
$\ssH\subset \ssG$ is cramped over the category of complex algebraic
representations if and only if $\comH\subset\comG$ is over the
category of complex smooth representations.

In \cite[Prop. 3.5.1]{WZ04}, Willenbring and Zuckerman prove the
following theorem:
\begin{thm}\label{WZ1}
  If $\ssG$ and $\ssH$ are semi-simple, and $\dim \ssH <\dim \ssG/P$
  for all parabolic subgroups $P$, then $\ssH$ is cramped.  In
  particular, if $\ssh\cong\fr{sl}_2\C$, then $\ssH$ is cramped if and
  only if $\ssh$ is not a summand of $\ssg$.
\end{thm}

In this paper, we will give necessary and sufficient conditions for an
inclusion of reductive groups to be cramped, using the geometry of
moment maps.  Peculiarly, it is not clear to the author how these
conditions imply Willenbring and Zuckerman's.  Instead, we use these
in conjunction with our results to prove a geometric result on the
intersection of adjoint orbits with Killing orthogonals to certain
subgroups.  An independent proof of this proposition (Proposition
\ref{WZ-cor}) would show that our results imply those of
Willenbring-Zuckerman.

\section{Moment Maps and Cramped Subgroups}

For all notation and basic notions of symplectic geometry, we refer the
reader to the book of Ana Cannas da Silva \cite{CdS}.

As usual, we consider $\comg^*$ as a Poisson manifold with the
Kostant-Kirillov Poisson structure: if $X,Y\in \fr g$, and we let
$i_X:\comg^*\to\R$ be the map given by pairing with $X$, then by
definition $\{i_X,i_Y\}=i_{[X,Y]}$.  Extending by the Leibniz rule and
continuity, we can take the Poisson bracket of any two functions. Note
that the linear functions of $\comg$ form a Lie subalgebra of the space of all
functions, which is isomorphic to $\comg$.

If $M$ is a symplectic (or more generally Poisson) manifold, a smooth
map \linebreak $\mu:M\to \comg^*$ is called {\bf a moment map} if it is Poisson
map.  In particular, if $\mu$ is a moment map, the functions
$\{\mu^*i_X\}_{X\in \comg}$ form a Lie subalgebra of functions under Poisson bracket.  After taking Hamiltonian vector fields, we obtain a map of Lie algebras from $\comg$ to vector fields on $M$. If $\comG$ is simply-connected, and $\mu$ proper, then this will
necessarily integrate to an action of $\comG$ on $M$.

If $N$ is any manifold with a $\comG$-action, then $M=T^*N$ has a canonical
symplectic structure, also equipped with a $\comG$-action.  This $\comG$-action
is Hamiltonian, that is, given by a moment map $\mu_N:M\to \comg^*$.  We
define $\mu_N$ as follows: for $n\in N,\xi\in T^*_nN$ and $X\in\comg$,
then $\langle\mu(n,\xi),X\rangle= \langle \xi,X_n\rangle$, where $X_n\in
T_nN$ is the vector given by differentiating the action of $G$ in the
direction of $X$.  Let $\im{N}$ be the image of $\mu_N$.

We can use moment maps in representation theory to gain information
using the philosophy of geometric quantization and geometric invariant
theory (for a readable heuristic discussion of this principle, see
\cite[\S 3.1]{GS}).

We will use a weak formulation 
of this connection:
\begin{thm} \label{KN}\emph{(Kirwan-Ness \cite{MFK})} Let
  $\comH\subset \comG$ be an inclusion of compact Lie groups.  Let
  $\la$ be a dominant integral weight of $\comG$ and $\mu$ a dominant
  integral weight of $\comH$, and let $\O_\la,\O_\mu$ the
  corresponding coadjoint orbits.  Let $\pi:\comg^*\to\comh^*$ be the
  natural projection given by restriction.  Then
\begin{equation*}
\Hom{W_{n\mu}}{V_{n\la}}{\comH}\neq \{0\}
\end{equation*}
for some integer $n$ if and only if $\O_\mu\subset \pi(\O_\la)$.
\end{thm}

This theorem will be the key to our understanding of $b(\ssg,\ssh)$.
We consider the $\comG$-space $N=\comG/\comH$.
\begin{thm} \label{main-theorem} The following are equivalent:
\begin{enumerate}
\item The subgroup $\ssH$ is cramped
\item The map $\mu_N:T^*N\to \comg^*$ is surjective.
\item For each maximal parabolic  subgroup $P\subset G$, there exists $x\in \comg$ s.t. 
$(x,\comh)=0$ and $\Stab_{\comG}x$ is conjugate to $P\cap K$.
\end{enumerate}
\end{thm}
Since $\im N=\mathrm{Ad}_{G}^*\cdot \fr h^{\perp}$, an inclusion of
reductive groups is cramped if and only if the corresponding inclusion
of algebras is.

Heuristically, we should consider $T^*N$ as a ``semi-classical limit''
of the $\comG$-representation $L^2(N)$, and $\im N$ as a
continuous analogue of the $\comG$-support of $L^2(N)$.  By Frobenius
reciprocity, that support is the same as the set of representations with
non-trivial $\comH$-invariants.  Thus, it is reasonable that $\mu_N$
being surjective reflects the fact that all representations are
``close'' to having $\comH$-invariants.

Let $X$ be a metric space with distance function $D:X\times X\to\R$.  For any subsets $A,B\subset X$, let
\begin{equation*}
  \delta(A,B)=\underset{b\in B}{\inf_{a\in A}} D(a,b).
\end{equation*}

\begin{proof}
  $(1)\Rightarrow(2)$: Assume there exists $\la\notin \im{N}$.  We can 
  assume $\la$ is dominant, possibly after applying conjugation by an 
  element of $\comG$.  Since $0$ not
  lying in $\pi(\O_\la)$ is an open condition on $\la$, we may
  furthermore assume that $m\la$ is integral for some $m$.  Of course,
  $\pi(\O_{m\la})$ also avoids $0$, so we may assume $\la$ is dominant
  and integral.  If $\Hom{W_\mu}{V_{m\la}}{\comH}\neq 0$ then by
  Theorem~\ref{KN} above, $\mu\in m\cdot\pi(\O_\la)$.  By homogeneity of
  $\im{N}$,
\begin{equation*}
  \del(m\la,\im{N})=\del(m\cdot \pi(\O_\la),0)=m\cdot\del(\pi(\comG\cdot
  \la),0). 
\end{equation*}
 
By the Weyl dimension formula, the set $\Hdim n=\{\mu\in\wlat_+^H|\dim
W_\mu < n\}$ is bounded for all integers $n$.  Let
$\eta_n=\sup_{\mu\in\Hdim n}\de(\mu,0)$.  Thus for each $n$, there
is another integer $N$ such that for $m> N$, the intersection $\Hdim
n\cap \left( m\cdot\pi(\O_\la)\right) $ is empty.  Thus, $V_{m\la}$ has no
non-trivial $H$-invariant subspaces of dimension less than $n$ and $b(\ssG,\ssH)=\infty$.

$(2)\Rightarrow(3)$: Consider the coadjoint orbit $\O$ through a point
with stabilizer $P$, that is, the dual of an orbit passing though the
center of the corresponding subalgebra $\fr p$.  Since $\mu_{N}$ is
surjective, the preimage of this orbit is non-empty and by the
transitivity of the action of $\comG$ on $N$, this preimage intersects
the fiber over $e\comH$, which is naturally identified with
$\comh^\perp$, the annihilator of $\comh$ in $\comg^*$ or
equivalently, under the Killing form on $\comg$.  Thus there is a
point in $\comh^\perp$ stabilized by a maximal parabolic subgroup
conjugate to $P$.

$(3)\Rightarrow(1)$: Since there is a unique orbit (up to scaling)
corresponding to each maximal parabolic subgroup, $\O_{\om_i}$ intersects
$\comh^\perp$.  Thus $0\in\pi(\O_{\om_i})$, and by Theorem~\ref{KN}, we
have $V_{m_i\om_i}^H\neq 0$ for some $m_i$.  The theorem now follows
from the arguments given by Willenbring and Zuckerman in
\cite[3.3.1,4.0.10]{WZ04}.  We sketch their argument for completeness:
if $x\in V_{m_i\om_i}^\ssH$, then multiplication by $x$ defines an
$H$-equivariant injection $V_{\la}\to V_{\la+m_i\om_i}$ by embedding
each representation in the coordinate ring of $\ssG/U$ (where $U$ is the
unipotent radical of $B$).  Thus, if $\la(\al_i)\geq m_i$ for any $i$,
then $V_\la$ has a $\ssH$-invariant subspace isomorphic to $V_{\la'}$
for some $\la'\leq \la$.  By induction,
\begin{equation*}
b(\ssG,\ssH)=\inf\{b(\la)|\la(\al_i^{\vee})<m_i\}<\infty.\qedhere
\end{equation*}
\end{proof}

Combining our results with those of Willenbring-Zuckerman, we obtain a seemingly new result on the geometry of Killing orthogonals to subgroups in semi-simple groups.

\begin{prop}\label{WZ-cor}
  Let $\comH\subset \comG$ be an inclusion of semi-simple, compact
  groups.  If every $\comG$-orbit $\mc O$ on $\comg$ satisfies $\dim
  \mc O > 2\dim\comH$, then for each orbit $\mc O$, we have $\mc
  O\cap\comh^{\perp}\neq \emptyset$, where $\comh^\perp$ is the
  orthogonal complement to $\comh$ under the Killing form.
\end{prop}
\begin{proof}
  By Theorem \cite[Prop. 3.5.8]{WZ04}, this inclusion
  $\comH\subset\comG$ is cramped since any adjoint orbit of minimal
  dimension is of the form $\comG/P\cap \comG$ for some maximal
  parabolic subgroup $P$ in the complexification $\ssG$ of $\comG$.  By
  Theorem~\ref{main-theorem}, the crampedness of the inclusion
  $\comH\subset\comG$ is equivalent to the desired result.
\end{proof}

If one could find an independent proof of this result, then it would
directly show that the necessary and sufficient conditions of
Theorem~\ref{main-theorem} imply Theorem \ref{WZ1}.  This seems to the
author to be more elegant, but we have found no proof of this result
notably simpler or shorter than the concatenation of Willenbring and
Zuckerman's with ours.

\section{Generalized Harish-Chandra modules}

Fix an inclusion of reductive Lie algebras $\tg\supset \ssg$.
Following Penkov and Zuckerman \cite{PZ04}, we call a $\tg$-module a
\textbf{generalized Harish-Chandra module} for the pair $(\tg,\ssg)$ if
as a $\ssg$-module if it is completely reducible with finite dimensional
isotypic components.  Obviously, all finite dimensional $\tg$-modules
lie in this category, but there are also many infinite dimensional ones
which are quite mysterious.

The original inspiration for considering cramped subgroups derives from
the theory of \ghc modules. Let $\ssh\subset\ssg$ is a reductive
subalgebra.  
\begin{prop}
  If $\ssh$ is cramped, then a \ghcm $V$ for $(\tg,\ssg)$ is also a \ghc module for
  $(\tg,\ssh)$ if and only if $V$ is finite-dimensional.
\end{prop}
\begin{proof}
  Obviously, $V$ is completely reducible as a $\ssh$-module, and all
  irreducible summands are finite-dimensional, so $V$ is a \ghc module if and
  only if all multiplicities of irreducible representations are finite.
  If $V$ is infinite dimensional, it has infinitely many isotypic
  components, and since $\ssh$ is cramped, each of these has an
  irreducible, non-trivial, $\ssh$-invariant subspace of dimension less
  than $b(\ssg,\ssh)$.  Since $\Hdim n$ is finite for all $n$, one of
  these representations has infinite multiplicity, and $V$ is not \ghc.
\end{proof}

Of course, the more general question arises, when is a $(\tg,\ssg)$ \ghcm
$V$ also a \ghc module for the pair $(\tg,\ssh)$ if $H$ is not cramped?

Using the techniques of symplectic geometry in the previous section, we will
give a characterization of \ghc modules for $(\tg,\ssg)$ which are also
\ghc for the pair $(\tg,\ssh)$ in terms of the image of the moment map
$\mu_N$.

Let $\supp_{\ssg}V\subset\wlat_+^{\ssg}$ be the $\ssg$-support of $V$, that
is, the set of weights $\la$ such that $\Hom{V_\la}{V}{\ssg}\neq\{0\}$.

For each $\gamma\in\R_{\geq 0}$, let 
\begin{math}
  \Mdis{\gamma}=\{\la\in \wlat^{\ssg}_+|\delta(\la,\im N)\leq\gamma \}
\end{math}.

\begin{thm}\label{supp-thm}
  A $(\tg,\ssg)$ \ghcm $V$, is a \ghcm for $(\tg,\ssh)$ if and only
  if for each $\gamma\in\R_{\geq 0}$, the intersection
  $\Mdis{\gamma}\cap \supp_{\ssg}V$ is finite. 
\end{thm}

We will first prove a lemma, which leads up to the proof of this theorem.  Let 
\begin{math}
  \Hdim n'=\{\la\in\wlat^\ssg_+|\supp_\ssH V_\la\cap\Hdim n\neq\emptyset\}.
\end{math}

\begin{lem}
  For each integer $n\in \Z$, and real number $\gamma$ there exists a
  real number $\beta\in\R_{\geq 0}$ and another integer $m\in \Z$, such
  that
  \begin{align*}
    \Hdim n'&\subset \Mdis\beta\\
    \Mdis\gamma&\subset\Hdim{m}'.
  \end{align*}
\end{lem}
\begin{proof}
  If $\la\in\Hdim n'$, then $\pi(\O_\la)\cap\Hdim n\neq \emptyset$.
  Applying $\pi$, we find that 
  \begin{equation*}
    \de(\la,\im N)=\de(\pi(\O_\la),0)\leq \eta_n
  \end{equation*}
  That is, $\Hdim n'\subset \Mdis{\eta_n}$.
  
  Note that $\Hdim 1'$ is a subsemigroup of full rank in
  $\im N\cap\wlat^{\ssg}_+$.  Thus, for each point in $\la\in
  \Mdis{\gamma}$ sufficiently far from the origin,
  there is a element $\la'\in\Hdim 1'$ such that
  $\la''=\la-\la'\in\Mdis{\gamma}$. By definition, $\la''$  is smaller
  than $\la$ in the
  standard poset structure of $\wlat^{\ssg}_+$.  Thus, for some
  $\eta_\gamma$, we have
  \begin{equation*}
    m_\gamma=\sup_{\la\in\Mdis{\gamma}}b(\la)=
    \underset{\de(\la,0)<\eta_\gamma}{\sup_{\la\in\Mdis{\gamma}}}b(\la)<\infty
  \end{equation*}
  since a ball of radius $\eta_\gamma$ is compact.  This is desired
  integer $m$.
\end{proof}

\begin{proof}[Proof of Theorem~\ref{supp-thm}]
  The module $V$ is \ghc if and only if the intersection
  $\Hdim n'\cap\supp_{\ssg}V$ is finite for each $n$.  Since each set
  $\Mdis\gamma$ contains and is contained in a set of the form $\Hdim
  n'$, these intersections are finite if and only if the intersections
  $\Mdis\gamma\cap\supp_{\ssg}V$ are also finite.
\end{proof}

Unfortunately, this is not a particularly easy criterion to apply, and
we aim to simplify it.  We can supply such a simplification, but at the
price of putting a restriction on our representations.  To do this, we
need the notion of asymptotic support and vagrancy.

Let $\asupp_{\ssg}V\subset\mathfrak{t}^*$ be the {\bf asymptotic
  support} of $V$, that is, the points $\la\in\mathfrak{t}^*$ such that
there is a sequence $a_i\in\R_{\geq 0}$ converging to 0, and a sequence
$\la_i\in\supp_{\ssg}V$ such that $\lim_{x\to\infty}a_i\la_i=\la$.

Define the {\bf vagrancy} of $V$ by
\begin{equation*}
      \irr(V) = \sup_{\la\in\supp_{\ssg}V}\delta(\la,\asupp_{\ssg}V)+\sup_{\la\in\asupp_{\ssg}V}\delta(\la,\supp_{\ssg}V).
\end{equation*}
\begin{cor}\label{asupp-cor}
  For all $(\tg,\ssg)$ \ghc modules $V$ such that $\irr (V)<\infty$, the
  following are equivalent:
  \begin{enumerate}
  \item $\im N\cap \asupp_{\ssg}V= \{0\}$.
  \item $V$ is \ghc for $(\tg,\ssh)$.
  \end{enumerate}
\end{cor}
Note that an essentially identical hypothesis, with similar conclusions
appears in the work of Kobayashi \cite{Kob03} in the much more difficult
context of arbitrary unitary representations of real groups.
\begin{proof}
  $(1)\Rightarrow(2)$: By the definition of $\irr (V)$, the set $\supp\,
  V$ is contained in the set of elements of distance less than or equal
  to $\irr (V)$ from $\asupp\, V$.  The intersection of this with
  $\Mdis{\gamma}$ for all $\gamma$ is compact (since it is closed, and
  bounded), and thus only has finitely many elements of $\supp\, V$.
  
  $(2)\Rightarrow(1)$: By the definition of $\irr (V)$, each ray in
  $\im N\cap\asupp\, V$ has infinitely many points of $\supp\, V$
  within distance $\irr (V)$ from it, and thus $\Mdis{\irr (V)}$ contains
  infinitely many elements of $\supp\,V$.
\end{proof}

As usual, this result raises more questions about the structure of
Harish-Chandra modules than it answers.  Foremost among them (in the
author's opinion) is the mysterious property of finite vagrancy.

\begin{question}
  Is $\irr (V)$ always finite?  Alternatively, is there a natural
  characterization of $V$ such that $\irr (V)$ is finite?
\end{question}

\section*{Acknowledgments}
\label{sec:acknowledgments}

I would like to thank Noah Snyder and Ivan Penkov for useful discussions.

\def\cprime{$'$}


\begin{thebibliography}{MFK94}

\bibitem[CdS01]{CdS}
Ana Cannas~da Silva.
\newblock {\em Lectures on symplectic geometry}, volume 1764 of {\em Lecture
  Notes in Mathematics}.
\newblock Springer-Verlag, Berlin, 2001.

\bibitem[GLS96]{GS}
Victor Guillemin, Eugene Lerman, and Shlomo Sternberg.
\newblock {\em Symplectic fibrations and multiplicity diagrams}.
\newblock Cambridge University Press, Cambridge, 1996.

\bibitem[MFK94]{MFK}
D.~Mumford, J.~Fogarty, and F.~Kirwan.
\newblock {\em Geometric invariant theory}, volume~34 of {\em Ergebnisse der
  Mathematik und ihrer Grenzgebiete (2) [Results in Mathematics and Related
  Areas (2)]}.
\newblock Springer-Verlag, Berlin, third edition, 1994.

\bibitem[Kob02]{Kob03}
Toshiyuki~Kobayashi. 
\newblock {\em Branching problems of unitary representations}, Proceedings of 
the ICM, Beijing, 2002, vol. 2, 615--627, Higher Ed. Press, Beijing, 2002. 

\bibitem[PZ04]{PZ04}
Ivan Penkov and Gregg Zuckerman.
\newblock Generalized {H}arish-{C}handra modules: a new direction in the
  structure theory of representations.
\newblock {\em Acta Appl. Math.}, 81(1-3):311--326, 2004.

\bibitem[WZ04]{WZ04}
Jeb~F. Willenbring and Gregg Zuckerman.
\newblock {Small semisimple subalgebras of semisimple Lie algebras}, 2004,
  arXiv:math.RT/0408302.

\end{thebibliography}
\end{document}